\documentclass[12pt, reqno]{amsart}

\usepackage{amsmath}
\usepackage{amsfonts}
\usepackage{amssymb}
\usepackage{amsthm}
\usepackage{eucal}
\usepackage[ansinew]{inputenc}
\usepackage{url}

\newtheorem*{ack}{Acknowledgements}
\numberwithin{equation}{section}

\begin{document}
\title {On a generalization of Jacobi's elegantissima}
\author[L.~Haine]{Luc Haine}
\address{Luc Haine, Institut de Recherche en Mathématique et Physique, Université catholique de
Louvain, Chemin du Cyclotron 2, 1348 Louvain-la-Neuve, Belgium}
\email{luc.haine@uclouvain.be}
\date{May 11, 2023; revised August 15, 2023}
\subjclass[2020]{Primary: 33E05; Secondary: 97I80}
\keywords{Elliptic functions, Landen's transformation}
\dedicatory{To the memory of Hermann Flaschka}
\begin{abstract}
We establish a generalization of Jacobi's elegantissima, which solves the pendulum equation. This amazing formula appears in lectures by the famous cosmologist Georges Lemaître, during the academic years 1955-1956 and 1956-1957. Our approach uses the full power of Jacobi's elliptic functions, in particular imagi-nary time is crucial for obtaining the result.
\end{abstract}
\maketitle
\tableofcontents

\section{Introduction}
In \emph{"Fundamenta Nova Theoriae Functionum Ellipticarum"}, Jacobi obtained the following formula for the "modular angle"
\begin{equation} \label{elegantissima}
\frac{\arcsin k}{4}=\arctan\sqrt{q}-\arctan\sqrt{q^3}+\arctan\sqrt{q^5}-\ldots,
\end{equation}
which is an inversion formula for
\begin{align}
K&=\int_0^1 \frac{dx}{\sqrt{(1-x^2)(1-k^2x^2)}},\; 0<k<1,\label{K}\\
K'&=\int_0^1 \frac{dx}{\sqrt{(1-x^2)(1-k'^2x^2)}},\;k'=\sqrt{1-k^2}\label{K'},
\end{align}
in terms of the so-called \emph{nome}
\begin{equation*}
q= e^{-\pi\frac{K'}{K}}.
\end{equation*}
Immediately following his formula (see \cite{J2} 40, formula (47)), Jacobi wrote \emph{"quae inter formulas elegantissima censeri debet"}. The only modern reference mentioning \eqref{elegantissima} seems to be \cite{WW}, Example 22.5.2.

During the two academic years 1955-1956 and 1956-1957, Georges Lemaître, the famous Belgian cosmologist, professor at the University of Louvain and father of the big bang theory, taught his \emph{"Leçons de mécanique. Le pendule"} \cite{L3}, \cite{L4}. For a beautiful account of Georges Lemaître's life and work, we refer the reader to \cite{L1}.

In the context of the oscillatory motion of a pendulum of length $l$ under the influence of the force of gravity, as explained in Section 2, Jacobi's formula can be seen as a formula giving the maximum deviation angle in terms of the period $4K\sqrt{l/g}$ of the oscillation and the \emph{imaginary period} $4K'\sqrt{l/g}$ of the \emph{complementary motion} obtained by inverting the direction of gravity, when the pendulum reaches the maximum angle. In his lectures \cite{L3}, \cite{L4}, Georges Lemaître derived a formula of a similar type, giving the angle of the oscillation at \emph{any} time, of which Jacobi's formula \eqref {elegantissima} becomes a particular case. According to his own words \emph{"La théorie que nous avons exposée peut être considérée comme une généralisation de l'élégantissime de Jacobi"}.

Lemaître's derivation of his generalization of Jacobi's elegantissima is based on a very ingenious infinite iteration of Landen's transformation \cite{L2}, a well known transformation of elliptic functions, see \cite{MM1}, \cite{MM2} for background. He started his lectures with elementary mechanical and geometric proofs of Jacobi's interpretation \cite{J1} of Poncelet theorem and Landen's transformation. Though not mentioned in his lectures, some of the material can also be found in Greenhill's book "The Applications of Elliptic Functions"; the French translation \cite{G2} we refer to, with a preface by Paul Appell, is a revised version of the original work. 

To the best of our knowledge, the infinite iteration of Landen's transformation to obtain the solution of the pendulum, does not seem to appear anywhere else in the literature. Though very intuitive, the passage to the limit is not carefully justified in \cite{L3}, \cite{L4}. The purpose of this note is to give a direct proof of Lemaître's formula, based on Jacobi's theory of elliptic functions, emphasizing \emph{the role played by ima-ginary time}. This is the content of Sections 3 and 4. In Section 5, we provide a sketch of Lemaître's proof. We don't follow his geometric arguments, using instead a standard version of Landen's transformation, to be found in \cite{WW}. Our hope is to attract the attention on Lemaître's elementary and geometric approach of elliptic functions.
\section{Jacobi's elliptic functions and the pendulum}
The motion of a pendulum of mass $m$ and length $l$ under the influence of the force of gravity is given by
\begin{equation} \label{pendule}
ml\ddot{\varphi}(t)=-mg\sin\varphi(t)\Leftrightarrow \ddot{\varphi}(t)=-\frac{\sin\varphi(t)}{R},\; R=\frac{l}{g},
\end{equation}
with $g$ the acceleration of the gravity, and $\varphi(t)$ the angle made with the descending vertical. Equation \eqref{pendule} has the following elementary solutions, the two equilibria given by
\begin{equation*}
\varphi(t)=0,\;\forall\;t\in\mathbb{R},\quad
\varphi(t)=\pi,\;\forall\;t\in\mathbb{R},
\end{equation*}
and, assuming $\varphi(0)=0$, the two doubly asymptotic solutions
\begin{equation} \label{das}
\varphi(t)=4\arctan e^{\pm\frac{t}{\sqrt{R}}}-\pi.
\end{equation}
All other solutions can be expressed in terms of Jacobi elliptic functions. Assuming $\varphi(0)=0$ and $\dot{\varphi}(0)>0$, as long as $\dot{\varphi}(t)>0$, one has
\begin{equation} \label{inversionp}
\sqrt{\frac{m}{2}}l \int_0^{\varphi(t)}\frac{d\varphi}{\sqrt{E-2mgl\sin^2\frac{\varphi}{2}}}=t,
\end{equation}
which follows from the conservation of energy
\begin{equation*}
\frac{ml^2\dot{\varphi}^2}{2}+2mgl\sin^2\frac{\varphi}{2}=E.
\end{equation*}

\emph{Oscillatory motions} occur when
\begin{equation*}
0<E<2mgl.
\end{equation*}
By the change of variable
\begin{equation*}
\sin\frac{\varphi}{2}=kx,
\end{equation*}
with
\begin{equation} \label{alpha}
0<k=\sin\frac{\alpha}{2}<1,
\end{equation}
where $0<\alpha<\pi$, gives the maximum angle of oscillation, i.e. $\dot{\varphi}(\alpha)=0$, equation \eqref{inversionp} becomes
\begin{equation*}
\int_0^{x(t)}\frac{dx}{\sqrt{(1-x^2)(1-k^2x^2)}}=\frac{t}{\sqrt{R}}.
\end{equation*}
Jacobi's elliptic sine function "sn" solves this inversion problem. Thus
\begin{equation} \label{so}
\sin\frac{\varphi(t)}{2}= k\;\mbox{sn}\Big(\frac{t}{\sqrt{R}},k\Big),\;
\cos\frac{\varphi(t)}{2}=\mbox{dn}\Big(\frac{t}{\sqrt{R}},k\Big),\;\forall\;t\in\mathbb{R},
\end{equation}
where the second formula follows from the identity $\mbox{dn}^2(z,k)=1-k^2\;\mbox{sn}^2(z,k)$, satisfied by Jacobi's elliptic sine "sn" and delta amplitude "dn" functions. Denoting by $T$ the period of the motion,
\begin{equation}\label{perosc}
\frac{T}{4\sqrt{R}}=\int_0^1\frac{dx}{\sqrt{(1-x^2)(1-k^2x^2)}}\Leftrightarrow T=4K\sqrt{R},
\end{equation}
with $K$ defined as in \eqref{K}.

For an \emph{oscillatory motion}, the \emph{complementary motion} is obtained by changing the sign of gravity, when the angle of oscillation reaches its maximum value $\alpha$. This motion follows the arc of circle complementary to the arc followed during the oscillation. It is characterized by
\begin{align*}
\ddot{\phi}(t)=-\frac{\sin\phi(t)}{R},\;
\phi(0)=\alpha-\pi,\;\dot{\phi}(0)=0.
\end{align*}
One immediately checks that the solution is given by
\begin{equation*}
\phi(t)=\varphi(it+K\sqrt{R})-\pi.
\end{equation*}
Since
\begin{equation*}
\sin^2 \frac{\phi(0)}{2}=\cos^2 \frac{\alpha}{2}=1-k^2,
\end{equation*}
from \eqref{perosc}, it follows that the period of the complementary motion is
\begin{equation*}
T'=4K'\sqrt{R},
\end{equation*}
with $K'$ as in \eqref{K'}. The mechanical interpretation of $T'$, using \emph{imaginary time}, is due to Paul Appell, see \cite{A}.

\emph{Circulatory motions} occur when
\begin{equation*}
E>2mgl\Leftrightarrow \frac{\dot{\varphi}^2(0)}{4}>\frac{1}{R}.
\end{equation*}
By the change of variable $x=\sin \frac{\varphi}{2}$,
equation \eqref{inversionp} becomes
\begin{equation*}
\int_0^{x(t)}\frac{dx}{\sqrt{(1-x^2)(1-k^2x^2)}}=\frac{t}{k\sqrt{R}},\quad 0< k=\frac{2}{\sqrt{R}\;\dot{\varphi}(0)}<1,
\end{equation*}
leading to
\begin{equation} \label{sc}
\sin \frac{\varphi(t)}{2}=\mbox{sn}\Big(\frac{t}{k\sqrt{R}},k\Big),\;\cos\frac{\varphi(t)}{2}=\mbox{cn}\Big(\frac{t}{k\sqrt{R}},k\Big),\; \forall\;t\in\mathbb{R},
\end{equation}
where the second formula follows from the identity $\mbox{cn}^2(z,k)+\mbox{sn}^2(z,k)=1$, satisfied by Jacobi's elliptic sine and cosine functions denoted "sn" and "cn". Using that $\mbox{sn}(z+2K,k)=-\mbox{sn}(z,k)$, one finds the period $T$ of the circulatory motion
\begin{equation} \label{pc}
T=2kK\sqrt{R}.
\end{equation}
\section{Lemaître's generalization of Jacobi's elegantissima}
For the oscillatory motion, defining
\begin{equation*}
\theta(t)=\pi+\varphi(t),
\end{equation*}
i.e. measuring the angle of oscillation from the ascending vertical, from \eqref{so} one gets
\begin{align}
\cos\frac{\theta(t)}{2}&=-\sin\frac{\varphi(t)}{2}=-k\;\mbox{sn}\Big(\frac{t}{\sqrt{R}},k\Big),\label{ocos}\\
\sin\frac{\theta(t)}{2}&=\cos\frac{\varphi(t)}{2}=\mbox{dn}\Big(\frac{t}{\sqrt{R}},k\Big)\label{osin}.
\end{align}
To derive from \eqref{ocos}, \eqref{osin} a nice formula for $e^{i\frac{\theta(t)}{2}}$, one is tempted to use the following infinite product formula (see \cite{WW} 22.5, for a proof)
\begin{multline} \label{ojt}
\mbox{dn}(z,k)+ik\;\mbox{sn}(z,k)=\\
\prod_{n=1}^\infty\frac{\big(1+(-1)^n q^{n-\frac{1}{2}}e^{-\frac{i\pi z}{2K}}\big)\big(1-(-1)^n q^{n-\frac{1}{2}}e^{\frac{i\pi z}{2K}}\big)}
{\big(1-(-1)^n q^{n-\frac{1}{2}}e^{-\frac{i\pi z}{2K}}\big)\big(1+(-1)^n q^{n-\frac{1}{2}}e^{\frac{i\pi z}{2K}}\big)},
\end{multline}
with
\begin{equation} \label{q}
q=e^{-\pi\frac{K'}{K}},
\end{equation}
then to take the logarithm to get $\theta(t)$. However, it is not immediately appropriate because of the imaginary exponentials which would appear in the infinite product. Thus we first rewrite the solution using Jacobi's imaginary transformation (see \cite{WW} 22.41) combined with standard formulas, which can all be obtained from the addition theorems for "sn", "cn" and "dn". Since
\begin{equation*}
\mbox{sn}(z-K,k)=-\frac{\mbox{cn}(z,k)}{\mbox{dn}(z,k)},\;
\mbox{dn}(iz+K',k')=\frac{k}{\mbox{dn}(iz,k')}=k\;\frac{\mbox{cn}(z,k)}{\mbox{dn}(z,k)},
\end{equation*}
we deduce that
\begin{equation} \label{pisn}
\mbox{sn}(z-K,k)=-\frac{1}{k}\;\mbox{dn}(iz+K',k').
\end{equation}
Similarly, since
\begin{equation*}
\mbox{dn}(z-K,k)=\frac{k'}{\mbox{dn}(z,k)},\;
\mbox{sn}(iz+K',k')=\frac{\mbox{cn}(iz,k')}{\mbox{dn}(iz,k')}=\frac{1}{\mbox{dn}(z,k)},
\end{equation*}
one has
\begin{equation}\label{pidn}
\mbox{dn}(z-K,k)=k'\;\mbox{sn}(iz+K',k').
\end{equation}
Combining \eqref{ocos}, \eqref{osin}, \eqref{pisn} and \eqref{pidn}, from \eqref{ojt}, with $k'$ instead of $k$, we obtain
\begin{align*}
e^{i\frac{\theta(t-K\sqrt{R})}{2}}&=\mbox{dn}\Big(i\;\frac{t}{\sqrt{R}}+K',k'\Big)+ik'\;\mbox{sn}\Big(i\;\frac{t}{\sqrt{R}}+K',k'\Big), \\
&=\prod_{n=1}^\infty\frac{\Big(1-i(-1)^n q'^{n-\frac{1}{2}}e^{\frac{\pi t}{2K'\sqrt{R}}}\Big)\Big(1-i(-1)^n q'^{n-\frac{1}{2}}e^{-\frac{\pi t}{2K'\sqrt{R}}}\Big)}
{\Big(1+i(-1)^n q'^{n-\frac{1}{2}}e^{\frac{\pi t}{2K'\sqrt{R}}}\Big)
\Big(1+i(-1)^n q'^{n-\frac{1}{2}}e^{-\frac{\pi t}{2K'\sqrt{R}}}\Big)},
\end{align*}
with
\begin{equation} \label{q'}
q'=e^{-\pi\frac{K}{K'}}.
\end{equation}
Since for $a\in\mathbb{R}$,
\begin{equation*}
e^{(-1)^n 2i \;arctan\;a}=\frac{1+(-1)^n i a}{1-(-1)^n i a},
\end{equation*}
taking the logarithm of the last formula, gives
\begin{multline}\label{legantissima}
\theta \Big(t-\frac{T}{4}\Big)=\\
\quad \sum_{n=1}^\infty(-1)^{n-1} 4 \Big\{\arctan \Big(q'^{n-\frac{1}{2}}\;e^{\frac{\pi t}{2K'\sqrt{R}}}\Big)+\arctan \Big(q'^{n-\frac{1}{2}}\;e^{-\frac{\pi t}{2K'\sqrt{R}}}\Big)\Big\},
\end{multline}
with $T=4K\sqrt{R}$ as in \eqref{perosc}, the period of the oscillatory motion. Evaluated at $t=0$, remembering that $\varphi(0)=0$,
\begin{equation*}
\theta\Big(-\frac{T}{4}\Big)=\pi+\varphi\Big(-\frac{T}{4}\Big)=\pi-\alpha,
\end{equation*}
with $\alpha$ the maximum angle of the oscillation, \eqref{legantissima} reduces to
\begin{equation*}
\pi-\alpha=8\;\sum_{n=1}^{\infty} (-1)^{n-1}\arctan\sqrt{q'^{2n-1}}.
\end{equation*}
Since the maximum angle of oscillation of the complementary motion is $\pi-\alpha$, permuting $K$ et $K'$ and remembering \eqref{alpha} and \eqref{q'}, we obtain
\begin{equation*}
2\arcsin k=\alpha=\pi-(\pi-\alpha)=8\sum_{n=1}^{\infty} (-1)^{n-1}\arctan\sqrt{q^{2n-1}},
\end{equation*}
with $q$ as in \eqref{q}. This is exactly \emph{Jacobi's elegantissima} \eqref{elegantissima}. Thus, formula \eqref{legantissima}, to be found in Lemaître \cite{L3}, \cite{L4}, deserves to be called \emph{the generalized elegantissima}.
\section{The circulatory motion}
Putting again
\begin{equation*}
\theta(t)=\pi+\varphi(t),
\end{equation*}
from \eqref{sc}, we now obtain
\begin{align}
\cos\frac{\theta(t)}{2}&=-\sin\frac{\varphi(t)}{2}=-\mbox{sn}\Big(\frac{t}{k\sqrt{R}},k\Big), \label{ccos}\\
\sin\frac{\theta(t)}{2}&=\cos\frac{\varphi(t)}{2}=\mbox{cn}\Big(\frac{t}{k\sqrt{R}},k\Big)\label{csin}.
\end{align}
Once again, we need to pass to imaginary time. Now we have
\begin{align*}
\mbox{cn}(z-K,k)&=k'\;\frac{\mbox{sn}(z,k)}{\mbox{dn}(z,k)},\\
\mbox{cn}(iz+K',k')&=-k\;\frac{\mbox{sn}(iz,k')}{\mbox{dn}(iz,k')}=-i k\;\frac{\mbox{sn}(z,k)}{\mbox{dn}(z,k)},
\end{align*}
and thus
\begin{equation} \label{picn}
\mbox{cn}(z-K,k)=i\frac{k'}{k}\;\mbox{cn}(iz+K',k').
\end{equation}
Combining \eqref{pisn}, \eqref{ccos}, \eqref{csin} and \eqref{picn}, we obtain
\begin{equation*}
e^{i\frac{\theta(t-kK\sqrt{R})}{2}}=\frac{\mbox{dn}\Big(i\frac{t}{k\sqrt{R}}+K',k'\Big)-k'\;\mbox{cn}\Big(i\frac{t}{k\sqrt{R}}+K',k'\Big)}{k}.
\end{equation*}
Therefore, we are led to use the following infinite product formula (see \cite{G2} 266, formula (39), but it can also be proven along the same lines as \eqref{ojt} following \cite{WW} 22.5, using an appropriate duplication formula)
\begin{multline*}
\frac{\mbox{dn}(z,k')-k'\;\mbox{cn}(z,k')}{k}=
\prod_{n=1}^\infty \frac{\big(1-q'^{n-\frac{1}{2}}e^{-\frac{i\pi z}{2K'}}\big)\big(1-q'^{n-\frac{1}{2}}e^{\frac{i\pi z}{2K'}}\big)}
{\big(1+q'^{n-\frac{1}{2}}e^{-\frac{i\pi z}{2K'}}\big)\big(1+q'^{n-\frac{1}{2}}e^{\frac{i\pi z}{2K'}}\big)},
\end{multline*}
with $q'=e^{-\pi \frac{K}{K'}}$. It gives
\begin{equation*}
e^{i\frac{\theta(t-kK\sqrt{R})}{2}}=
\prod_{n=1}^\infty \frac{\Big(1+iq'^{n-\frac{1}{2}}e^{\frac{\pi t}{2kK'\sqrt{R}}}\Big)\Big(1-iq'^{n-\frac{1}{2}}e^{-\frac{\pi t}{2kK'\sqrt{R}}}\Big)}
{\Big(1-iq'^{n-\frac{1}{2}}e^{\frac{\pi t}{2kK'\sqrt{R}}}\Big)\Big(1+iq'^{n-\frac{1}{2}}e^{-\frac{\pi t}{2kK'\sqrt{R}}}\Big)}.
\end{equation*}
Taking the logarithm, we now obtain the following formula
\begin{multline}
\theta\Big(t-\frac{T}{2}\Big)=\\
\sum_{n=1}^\infty 4\Big\{\arctan\Big(q'^{n-\frac{1}{2}}\;e^{\frac{\pi t}{2kK'\sqrt{R}}}\Big)-\arctan\Big(q'^{n-\frac{1}{2}}\;e^{-\frac{\pi t}{2kK'\sqrt{R}}}\Big)\Big\} \label{lcirc},
\end{multline}
with  $T=2kK\sqrt{R}$ as in \eqref{pc}, the period of the circulatory motion.
\section{A sketch of Lemaître's proof}
In this Section, we don't follow Lemaître's arguments in detail, for which we refer the reader to \cite{L4}, the second version of his lectures during the academic year 1956-1957. Our aim is to sketch the spirit of his approach, but we use freely a standard form of Landen's transformation to be found in \cite{WW}. We present his arguments in four steps. At each step, we refer to the sections in the table of contents of \cite{L4}, where elementary proofs using plane geometry are given.\\

\underline{Step 1}. The basic idea (see \cite{L4} 10, 11) is to use Landen's transformation to express the circulatory motion as a sum of two \emph{new} circulatory motions \emph{identical up to a shift by half a period}. Precisely, if $\theta(t)$ is a circulatory motion with period $T=2kK\sqrt{R}$ as in \eqref{pc}, there exists a circulatory motion $\theta_1(t)$ with period $T_1=2k_1K_1\sqrt{R_1}$ such that
\begin{equation} \label{dcirc}
\theta\Big(\frac{R}{R_1}t\Big)=\theta_1(t)+\theta_1 \Big(t-\frac{T_1}{2}\Big),
\end{equation}
and the period $T_1$ of $\theta_1(t)$ is related to the period $T$ of $\theta(t)$ as follows
\begin{equation} \label{T1}
T_1=2\frac{R_1}{R}T.
\end{equation}
To establish this result, using \eqref{ccos}, \eqref{csin} and the addition theorem for "sn" and "cn", we rewrite \eqref{dcirc} as follows
\begin{align} \label{dtcirc}
\cos\frac{\theta\Big(\frac{R}{R_1}t\Big)}{2}&=
\cos \frac{\theta_1(t)}{2}\cos \frac{\theta_1(t-k_1K_1\sqrt{R_1})}{2}\nonumber\\
&\quad -\sin\frac{\theta_1(t)}{2}\sin \frac{\theta_1(t-k_1K_1\sqrt{R_1})}{2}, \nonumber\\
&=\mbox{sn}\Big(\frac{t}{k_1\sqrt{R_1}},k_1\Big)\mbox{sn}\Big(\frac{t}{k_1\sqrt{R_1}}-K_1,k_1\Big)\nonumber\\
&\quad -\mbox{cn}\Big(\frac{t}{k_1\sqrt{R_1}},k_1\Big)\mbox{cn}\Big(\frac{t}{k_1\sqrt{R_1}}-K_1,k_1\Big),\nonumber\\
&=-(1+k'_1)\frac{\mbox{sn}\Big(\frac{t}{k_1\sqrt{R_1}},k_1\Big)\mbox{cn}\Big(\frac{t}{k_1\sqrt{R_1}},k_1\Big)}{\mbox{dn}\Big(\frac{t}{k_1\sqrt{R_1}},k_1\Big)}.
\end{align}
Now Landen's transformation (see \cite{WW} 22.42) tells us that \eqref{dtcirc} can be written as follows
\begin{equation} \label{dlcirc}
\cos\frac{\theta\Big(\frac{R}{R_1}t\Big)}{2}=-\mbox{sn}\Big(\frac{(1+k'_1)t}{k_1\sqrt{R_1}},k\Big),
\end{equation}
with
\begin{equation}\label{k1}
k=\frac{1-k'_1}{1+k'_1}=\Big(\frac{k_1}{1+k'_1}\Big)^2.
\end{equation}
Thus assuming, as prescribed by \eqref{ccos}, that
\begin{equation*}
\cos\frac{\theta(t)}{2}=-\mbox{sn}\Big(\frac{t}{k\sqrt{R}},k\Big),
\end{equation*}
to satisfy \eqref{dlcirc}, we must pick
\begin{equation} \label{R_1}
R_1=\frac{R}{k},
\end{equation}
and choose $k_1$ as function of $k$ as in \eqref{k1}. It is a standard result (see \cite{WW} 22.42) that
\begin{equation} \label{KK1}
K=\frac{(1+k'_1)K_1}{2},
\end{equation}
with $K_1$ as in \eqref{K} by replacing $k$ with $k_1$, hence
\begin{equation*}
T_1=2k_1 K_1\sqrt{R_1}=4K\sqrt{R}=2\frac{R_1}{R}T,
\end{equation*}
which establishes \eqref{T1}.

Iterating $j$ times the previous construction, we obtain a circulatory motion $\theta_j(t)$ with period
\begin{equation} \label{Tj}
T_j=2k_jK_j\sqrt{R_j}=2^j\frac{R_j}{R}T,
\end{equation}
such that the circulatory motion $\theta(t)$ with period $T=2kK\sqrt{R}$ we started with, can be expressed as follows
\begin{equation}\label{dfc}
\theta (t)=\sum_{n=-2^{j-1}}^{2^{j-1}-1} \theta_j\Big(\frac{R_j}{R}(t+n T)\Big),
\end{equation}
with $R_j,k_j$ inductively defined by
\begin{equation} \label{ifc}
R_j=\frac{R_{j-1}}{k_{j-1}},\; k'_j=\sqrt{1-k_j^2}=\frac{1-k_{j-1}}{1+k_{j-1}}, j\geq 1,
\end{equation}
with $R_0=R, k_0=k$.

Lemaître also observes that connecting the $2^j$ successive circulatory motions in \eqref{dfc} by segments, one gets a closed Poncelet polygon (see \cite{L4} 9, 15, 22). The next step is to establish that $\lim_{j\to\infty} R_j=R_\infty<\infty$.\\

\underline{Step 2}. We follow \cite{L4} (13, 14, 15).  From \eqref{ifc}, one easily establishes that
\begin{equation*}
\frac{R_n}{R_{n+1}}=\frac{\sqrt{R_{n-1}R_n}}{\frac{R_{n-1}+R_{n}}{2}}, \;n\geq 1.
\end{equation*}
Putting
\begin{equation} \label{a0b0}
b_0=R <a_0=R_1,
\end{equation}
as shown by Gauss \cite{G1}, the two sequences
\begin{equation*}
a_{n}=\frac{a_{n-1}+b_{n-1}}{2}, \;b_n=\sqrt{a_{n-1}b_{n-1}}, \;n\geq 1,
\end{equation*}
converge to a common limit denoted by $M(a_0,b_0)$, the so-called arithmetic-geometric mean
\begin{equation*}
M(a_0,b_0)=\lim_{n\to\infty}a_n=\lim_{n\to\infty}b_n,
\end{equation*}
and moreover
\begin{equation} \label{gauss}
\frac{\pi}{2M(a_0,b_0)}=\int_0^{\frac{\pi}{2}}\frac{d\phi}{\sqrt{a_0^2\cos^2\phi+b_0^2\sin^2\phi}}.
\end{equation}
We refer the reader to \cite{MM2} 2.3, for a beautiful account of the arithmetic-geometric mean, and to \cite{BB} and \cite{MM1} for more historical background and updated research on the subject.

We have
\begin{equation*}
\frac{R_1}{R_2}=\frac{\sqrt{RR_1}}{\frac{R+R_1}{2}}=\frac{b_1}{a_1},
\end{equation*}
and, by induction
\begin{equation*}
\frac{R_n}{R_{n+1}}=\frac{\sqrt{\frac{R_{n-1}}{R_n}}}{\frac{1}{2}(1+\frac{R_{n-1}}{R_n})}=\frac{\sqrt{a_{n-1}b_{n-1}}}{\frac{a_{n-1}+b_{n-1}}{2}}=\frac{b_n}{a_n}.
\end{equation*}
Hence, using that
\begin{equation*}
b_n^2b_{n+1}^2\ldots b_{n+j}^2=a_{n-1}a_n\ldots a_{n+j-1}b_{n-1}b_n\ldots b_{n+j-1},
\end{equation*}
we get
\begin{align*}
\frac{R_n}{R_{n+j}}&=\frac{R_n}{R_{n+1}}\frac{R_{n+1}}{R_{n+2}}\ldots\frac{R_{n+j-1}}{R_{n+j}}
=\frac{b_nb_{n+1}\ldots b_{n+j-1}}{a_na_{n+1}\ldots a_{n+j-1}}\\
&=\frac{a_{n-1}b_{n-1}}{b_{n+j}^2}
=\frac{b_n^2}{b_{n+j}^2}.
\end{align*}
In particular, for $n=0$, we have
\begin{equation*}
\frac{R}{R_j}=\frac{b_0^2}{b_j^2},
\end{equation*}
which gives
\begin{equation*}
R_\infty=\lim_{j\to\infty}R_j=\frac{R}{b_0^2}\lim_{j\to \infty} b_j^2=R\Big(\frac{M(a_0,b_0)}{b_0}\Big)^2.
\end{equation*}
Remembering \eqref{R_1} and \eqref{a0b0}, $\frac{b_0}{a_0}=\frac{R}{R_1}=k$, we compute from \eqref{gauss} that
\begin{equation*}
\frac{\pi}{2M(a_0,b_0)}=\frac{k}{R}\int_0^{\frac{\pi}{2}}\frac{d\phi}{\sqrt{1-k'^2\sin^2\phi}}=\frac{k}{R}K',
\end{equation*}
with $K'$ as in \eqref{K'}, leading to
\begin{equation} \label{rl}
R_\infty=\Big(\frac{\pi}{2kK'}\Big)^2R.
\end{equation}

The next step, the passage to the limit in \eqref{dfc}, is not carefully justified in Lemaître's syllabi \cite{L3}, \cite{L4}.  We won't attempt to justify it, since we have obtained a different proof of his result in Section 4.\\

\underline{Step 3}. This is developed in \cite{L4} (16, 17). Remembering \eqref{Tj}, during the iteration, the period essentially doubles at each step. Since the limit radius $R_\infty$ is finite, at the limit we obtain the doubly asymptotic motion \eqref{das} of the pendulum
\begin{equation*}
\lim_{j\to\infty} \theta_j(t)=\theta_\infty(t)=4\arctan e^{\frac{t}{\sqrt{R_\infty}}},
\end{equation*}
with $\theta_\infty(0)=\pi,\dot\theta_\infty(0)>0$. Taking \emph{formally} the limit of \eqref{dfc}, we get
\begin{equation*}
\theta(t)=\sum_{n=-\infty}^\infty \theta_\infty \Big(\frac{R_\infty}{R}(t+nT)\Big),
\end{equation*}
i.e., using \eqref{rl}, we obtain
\begin{equation} \label{fl}
\theta(t)=\sum_{n=-\infty}^\infty 4\arctan e^{\frac{\pi}{2kK'\sqrt{R}}(t+n T)}.
\end{equation}
However, this series is not convergent. Putting
\begin{equation*}
q'=e^{-\pi\frac{K}{K'}},
\end{equation*}
and remembering \eqref{pc}, we can rewrite \eqref{fl} as follows
\begin{align*}
&\theta\Big(t-\frac{T}{2}\Big)=\sum_{n=1}^\infty 4\arctan \Big(q'^{n-\frac{1}{2}}\;e^{\frac{\pi t}{2kK'\sqrt{R}}}\Big)\nonumber\\
&+\sum_{n=1}^\infty 4\arctan \Big(q'^{-(n-\frac{1}{2})}\;e^{\frac{\pi t}{2kK'\sqrt{R}}}\Big),\nonumber\\
\quad &=\sum_{n=1}^\infty 4\Big\{\arctan\Big(q'^{n-\frac{1}{2}}\;e^{\frac{\pi t}{2kK'\sqrt{R}}}\Big)-\arctan\Big(q'^{n-\frac{1}{2}}\;e^{-\frac{\pi t}{2kK'\sqrt{R}}}\Big)\Big\} \label{lcirc},
\end{align*}
where each term in the second sum is recomputed modulo $2\pi$ using
\begin{equation*}
4\arctan\frac{1}{x}=2\pi-4\arctan x,
\end{equation*}
which agrees with formula \eqref{lcirc} established in Section 4. Finally, and it is the last step, Lemaître obtains his generalization of Jacobi's elegantissima by another application of Landen's transformation.\\

\underline{Step 4}. Following \cite{L4} (12, 18), we consider \emph{two identical circulatory motions} $\theta_1(t)$, with period $T_1=2k_1K_1\sqrt{R_1}$, up to a shift by half a period, going now in \emph{opposite} directions, i.e. we look at
\begin{equation*}
\theta(t)=\theta_1(t)-\theta_1(t-k_1K_1\sqrt{R_1}).
\end{equation*}
By a similar computation as in \eqref{dtcirc} and \eqref{dlcirc}, using \eqref{ccos} and \eqref{csin},  we get
\begin{equation*}
\cos\frac{\theta(t)}{2}=-\frac{1-k'_1}{1+k'_1}\;\mbox{sn}\Big(\frac{(1+k'_1)t}{k_1\sqrt{R_1}}, \frac{1-k'_1}{1+k'_1}\Big),
\end{equation*}
i.e. defining
\begin{equation} \label{kR}
k=\frac{1-k'_1}{1+k'_1}=\Big(\frac{k_1}{1+k'_1}\Big)^2 \quad \mbox{and} \quad R=k R_1,
\end{equation}
we obtain
\begin{equation*}
\cos\frac{\theta(t)}{2}=-k\;\mbox{sn}\Big(\frac{t}{\sqrt{R}},k\Big).
\end{equation*}
Remembering \eqref{ocos}, this shows that $\theta(t)$ is now an oscillatory motion with period $T=4K\sqrt{R}$ as in \eqref{perosc}. Using \eqref{KK1} and \eqref{kR}, one finds
\begin{equation} \label{T1T}
T_1=2k_1K_1\sqrt{R_1}=\frac{k_1}{1+k'_1}4K\sqrt{R_1}=4K\sqrt{R}=T,
\end{equation}
i.e. the oscillatory motion $\theta(t)$ and the circulatory motion $\theta_1(t)$ have now the same period. It is a standard result  (see \cite{WW}, Example 22.4.5) that
\begin{equation} \label{K'K'1}
K'=(1+k'_1)K'_1,
\end{equation}
with $K'_1$ as in \eqref{K'} by replacing $k'$ with $k'_1$. Hence, using \eqref{kR}, we deduce
\begin{equation} \label{R1R}
2k_1K'_1\sqrt{R_1}=\frac{k_1}{1+k'_1}2K'\sqrt{R_1}=2K'\sqrt{R}.
\end{equation}
Using \eqref{KK1} and \eqref{K'K'1}, we also have
\begin{equation}\label{pq1q}
q'_1=e^{-\pi\frac{K_1}{K'_1}}=e^{-2\pi\frac{K}{K'}}=(q')^2.
\end{equation}
From \eqref{lcirc}, using \eqref{T1T}, \eqref{R1R} and \eqref{pq1q}, we obtain
\begin{multline*}
\theta_1\Big(t-\frac{T_1}{2}\Big)=\\
\sum_{n=1}^\infty 4\Big\{\arctan \Big(q'^{2n-1}\;e^\frac{\pi t}{2K'\sqrt{R}}\Big)-\arctan \Big(q'^{2n-1}\;e^{-\frac{\pi t}{2K'\sqrt{R}}}\Big) \Big\},
\end{multline*}
hence, by an easy computation, we get
\begin{align*}
&\theta\Big(t-\frac{T}{4}\Big)=\theta_1\Big(t-\frac{T_1}{2}+\frac{T_1}{4}\Big)-\theta_1\Big(t-\frac{T_1}{2}-\frac{T_1}{4}\Big)\\
&=\sum_{n=1}^\infty(-1)^{n-1} 4\Big\{\arctan \Big(q'^{n-\frac{1}{2}}\;e^{\frac{\pi t}{2K'\sqrt{R}}}\Big)+\arctan \Big(q'^{n-\frac{1}{2}}\;e^{-\frac{\pi t}{2K'\sqrt{R}}}\Big)\Big\},
\end{align*}
which is \emph{the generalized elegantissima} as in \eqref{legantissima}.

\begin{ack}
\emph{We thank the referees for helpful comments and for pointing out the relevance of references \cite{BB} and \cite{MM1} to this work.}
\end{ack}


\begin{thebibliography}{xx} 
\bibitem{A} P. Appell, \emph{Sur une interprétation des valeurs imaginaires du temps en Mécanique}, Comptes Rendus Hebdomadaires des Séances de l'Académie des Sciences, t.87 (30 déc. 1878) 1074-1077.
\bibitem{BB} J.M. Borwein and P.B. Borwein, \emph{Pi and the AGM: A Study in Analytic Number Theory and Computational Complexity}, Canadian Mathematical Society Series of Monographs and Advanced Texts, John Wiley \& Sons Inc., New York (1987). 
\bibitem{G1} C.F. Gauss, \emph{Arithmetisch Geometrisches Mittel}, in Werke, vol.3 (1799) 361-432.
\bibitem{G2} A.G. Greenhill, \emph{Les fonctions elliptiques et leurs applications}, Editions Jacques Gabay (1895).
\bibitem{J1} C.G.J. Jacobi, \emph{Ueber die anwendung der elliptischen transcendenten auf ein bekanntes problem der elementargeometrie}, Crelle Journal für die reine und angewandte Mathematik 3 (1828) 376-389.
\bibitem{J2} C.G.J. Jacobi, \emph{Fundamenta Nova Theoriae Functionum Ellipticarum}, Bornträger, Königsberg (1829).
\bibitem{L1} D. Lambert, \emph{Un atome d'univers}, La vie et l'oeuvre de Georges Lemaître, Editions Lessius, Bruxelles et Editions Racines, Bruxelles (2000).
\bibitem{L2} J. Landen, \emph{An investigation of a general theorem for finding the length of any arc of any conic hyperbola, by means of two elliptic arcs, with some other new and useful theorems deduced therefrom}, Philos. Trans. Royal Soc. London 65 (1775) 283-289.
\bibitem{L3} G. Lemaître, \emph{Leçons de mécanique. Le pendule}. Notes de cours (1955-1956)
\url{https://archives.uclouvain.be/ark:/33176/dli000000aGtr#?c=0&m=0&s=0&cv=67&xywh=-479%2C-1%2C5309%2C2650}.
\bibitem{L4} G. Lemaître, \emph{Leçons de mécanique. Le pendule}. Notes de cours (1956-1957)
\url{https://archives.uclouvain.be/ark:/33176/dli000000aGNf#?c=0&m=0&s=0&cv=2&xywh=-605%2C0%2C5279%2C2634}
\bibitem{MM1} D.V. Manna and V.H. Moll, \emph{Landen survey}, Probability, Geometry and Integrable Systems, MSRI Publications, For Henry McKean's Seventy-Fifth Birthday (Edited by M. Pinsky and B. Birnir) vol. 55 (2007) 287-319.
\bibitem{MM2} H. McKean and V. Moll, \emph{Elliptic curves}, Function theory, geometry, arithmetic, Cambrige University Press (1997).
\bibitem{WW} E.T. Whittaker and G.N. Watson, \emph{A course of modern analysis}, Fifth edition edited and prepared for publication by Victor H. Moll, Cambridge University Press (2021).
\end{thebibliography}
\end{document}